\documentclass[10pt]{article}
\usepackage{amsmath}
\usepackage{amsfonts}
\usepackage{amssymb}
\usepackage{comment}
\long\def\proof#1{\removelastskip\vskip\baselineskip\relax\noindent{\it
Proof\if!#1!\else\ \ignorespaces#1\fi.\ }\ignorespaces}

\newcommand{\lgs}[2]{\mbox{$\left(\frac{#1}{#2}\right)$}}

\renewcommand{\th}{\theta}

\newcommand{\G}{\Gamma}

\newcommand{\eps}{\varepsilon}

\newcommand{\Proof}{{\it Proof. \/}}
\newcommand{\squareforqed}{\hbox{\rlap{$\sqcap$}$\sqcup$}}
\newcommand{\qed}{\ifmmode\squareforqed\else{\unskip\nobreak\hfil
\penalty50\hskip1em\null\nobreak\hfil\squareforqed
\parfillskip=0pt\finalhyphendemerits=0\endgraf}\fi}

\newcommand{\fp}{\qed\removelastskip\vskip\baselineskip\relax}

\newtheorem{theorem}{Theorem}[section]

\newtheorem{proposition}[theorem]{Proposition}

\newtheorem{definition}[theorem]{Definition}

\newcommand{\litem}{\par\noindent\dimen0=\parindent%
    \advance\dimen0 by-4pt
               \hangindent=\dimen0\ltextindent}

\newcommand{\ltextindent}[1]{\hbox to \hangindent{#1\hss}\ignorespaces}
\newcommand{\ltextjndent}[1]{\hbox to \hangindent{#1\hss}\ignorespaces\kern-1ex}

\renewcommand{\pmod}[1]{\allowbreak\ ({\rm{mod}}\,\,#1)}

\begin{document}
\pagestyle{plain}

\title{Factorizations of Eisenstein Series of Level up to $4$}
\author{Henri Cohen,\\
Universit\'e de Bordeaux,\\
LFANT, IMB, U.M.R. 5251 du C.N.R.S,\\
351 Cours de la Lib\'eration,\\
33405 TALENCE Cedex, FRANCE}

\maketitle
\begin{abstract}
  We show that the most standard Eisenstein series such as $E_4(\tau)$
  or $2E_2(2\tau)-E_2(\tau)$, and also the function $\th^2(\tau)$, are in a
  natural way the product of two conjugate Eisenstein series of half their
  weight and double their level, as well as a number of similar elementary
  identities for $E_6$ and Eisenstein series of levels $2$, $3$, and $4$.
\end{abstract}

\smallskip

\section{Introduction and Notation}

The goal of this paper is to show that the most basic Eisenstein series
of level up to $4$ have natural factorizations into a product of two
conjugate Eisenstein-type series of half their weight and double their
level as well as similar identities, and to study the corresponding functions.
No proofs are explicitly given since in all cases it is simply a matter of
working in a suitable finite-dimensional space of modular forms on some
congruence subgroup of the modular group and checking a few initial Fourier
coefficients.

The motivation for this paper comes from the fact that the functions
$F_{N,j}$ that we will introduce are essentially normalized Hauptmoduln
for several standard non-cocompact arithmetic triangle groups, see
\cite{Coh}.

We use the standard notation of modular forms and functions: we will use
in particular the functions $E_2$, $E_4$, $E_6$, $\eta$, and $\th$,
and we will also need the Eisenstein series of weight $1$ and levels $3$ and
$4$:

$$E_{1,-3}(\tau)=1/6+\sum_{n\ge1}\bigl(\sum_{d\mid n}\lgs{-3}{d}\bigr)q^n\text{\quad and\quad}
E_{1,-4}(\tau)=1/4+\sum_{n\ge1}\bigl(\sum_{d\mid n}\lgs{-4}{d}\bigr)q^n\;,$$
and we recall that $\th^2=4E_{1,-4}$.
  
For $N=2$, $3$, and $4$ we also define the forms
\begin{align*}E_{N,2}(\tau)&=(NE_2(N\tau)-E_2(\tau))/(N-1)\in M_2(\G_0^*(N))\text{\quad and}\\
  E_{N,4}(\tau)&=(N^2E_4(N\tau)-E_4(\tau))/(N^2-1)\in M_4(\G_0^*(N))\;,\end{align*}
\noindent  
where as usual $\G_0^*(N)$ denotes the group generated by $\G_0(N)$
and the Fricke involution $\tau\mapsto-1/(N\tau)$, but we will not use
the function $E_{4,4}$. In levels $3$ and $4$ we will need additional
modular forms which will be introduced when needed.

We could present our results in a more or less uniform manner, but
it is more readable to give them level by level, referring to the above
for definitions of the functions which are used.

\section{Results in Level $1$}

The basic formula in level $1$ is
$$E_4^3-E_6^2=1728\eta^{24}\;.$$
This can be factored in two ways and leads to the following:
\begin{definition}
  We define
    $$F_{1,j}(\tau)=E_6(\tau)-(-1)^j\sqrt{-1728}\eta^{12}(\tau)\text{\quad and\quad}G_{1,j}(\tau)=E_4(\tau)-12\rho^j\eta^8(\tau)\;,$$
    where as usual $\rho=(-1+\sqrt{-3})/2$.\end{definition}

\begin{proposition}
  We have
  \begin{align*}
    &F_{1,j}(\tau+1)=F_{1,j+1}(\tau)\;,\quad F_{1,j}(-1/\tau)=\tau^6F_{1,j+1}(\tau)\;,\\
    &F_{1,j}(\tau)F_{1,j+1}(\tau)=E_4^3(\tau)\;,\\
    &G_{1,j}(\tau+1)=G_{1,j+1}(\tau)\;,\quad G_{1,j}(-1/\tau)=\tau^4G_{1,j}(\tau)\;,\\
    &G_{1,j}(\tau)G_{1,j+1}(\tau)G_{1,j+2}(\tau)=E_6^2(\tau)\;.\end{align*}
\end{proposition}

\Proof Trivial.\fp

\begin{theorem}[Level 1]
  \begin{enumerate}
  \item The functions
    \begin{align*}f_{1,j}(\tau)=2E_2(\tau)&-((3+(-1)^j\sqrt{-3})/6)E_2(\tau/2)\\
      &-((3-(-1)^j\sqrt{-3})/6)E_2((\tau+1)/2)\end{align*}
  have the following properties:
  \begin{align*}
   &f_{1,j}^3(\tau)=F_{1,j}(\tau)\;,\\
   &f_{1,j}(\tau)f_{1,j+1}(\tau)=E_4(\tau)\;,\quad f_{1,j}(((-1)^j+\sqrt{-3})/2)=0\;,\\
   &f_{1,j}(\tau+1)=f_{1,j+1}(\tau)\;,\quad f_{1,j}(-1/\tau)=((-1-(-1)^j\sqrt{-3})/2)\tau^2f_{1,j+1}(\tau)\;.\end{align*}
\item The functions
  $$g_{1,j}(\tau)=(3/2)E_2(\tau)-(1/4)\sum_{\substack{0\le m\le 2\\m\not\equiv j\pmod{3}}}E_2((\tau+m)/3)$$
  have the following properties:
  \begin{align*}
    &g_{1,j}^2(\tau)=G_{1,j}(\tau)\;,\\
    &g_{1,j}(\tau)g_{1,j+1}(\tau)g_{1,j+2}(\tau)=E_6(\tau)\;,\quad g_{1,j}(-j+\sqrt{-1})=0\;,\\
    &g_{1,j}(\tau+1)=g_{1,j+1}(\tau)\;,\quad g_{1,j}(-1/\tau)=\eps_j\tau^2g_{1,j}(\tau)\;,\end{align*}
  with $\eps_j=1$ if $3\mid j$ and $\eps_j=-1$ if $3\nmid j$.
  \end{enumerate}
  \end{theorem}

\Proof All the functions involved are modular over $\G(2)$ or $\G(3)$, so the
proof is a simple verification of a few coefficients. Evidently this is not an
explanation of the theorem, especially since a similar theorem is valid in
levels $2$, $3$, and $4$.\fp

In particular this shows that $E_4$ is the product of two conjugate Eisenstein
series of weight $2$ and level $2$, and that $E_6$ is the product of three
conjugate Eisenstein series of weight $2$ and level $3$.

\section{Results in Level $2$}  

For simplicity, we set $\eta_2(\tau)=\eta(\tau)\eta(2\tau)$, of weight $1$
(not to be confused with $\eta(2\tau)$ itself).

In this level, we have the following eta quotient identities:
$$E_{2,2}^2(\tau)-E_{2,4}(\tau)=128\eta^{16}(2\tau)/\eta^8(\tau)\;,\quad
E_{2,2}^2(\tau)+E_{2,4}(\tau)=2\eta^{16}(\tau)/\eta^8(2\tau)\;.$$
The basic formula in level $2$ is
$$E_{2,2}^4-E_{2,4}^2=256\eta_2^8\;.$$
This can be factored in two ways and leads to the following:
\begin{definition}
  We define
 $$F_{2,j}(\tau)=E_{2,4}(\tau)-(-1)^j\sqrt{-256}\eta_2^4(\tau)\text{\quad and\quad}G_{2,j}(\tau)=E_{2,2}(\tau)-4i^j\eta_2^2(\tau)\;,$$
 where $\eta_2$ has been defined above.
\end{definition}

\begin{proposition}
  We have
  \begin{align*}
    &F_{2,j}(\tau+1)=F_{2,j+1}(\tau)\;,\quad F_{2,j}(-1/(2\tau))=-4\tau^4F_{2,j+1}(\tau)\;,\\
    &F_{2,j}(\tau)F_{2,j+1}(\tau)=E_{2,2}^4(\tau)\;,\\
    &G_{2,j}(\tau+1)=G_{2,j+1}(\tau)\;,\quad G_{2,j}(-1/(2\tau))=-2\tau^2G_{2,j}(\tau)\;,\\
    &G_{2,j}(\tau)G_{2,j+1}(\tau)G_{2,j+2}(\tau)G_{2,j+3}(\tau)=E_{2,4}^2(\tau)\;.\end{align*}\end{proposition}

\Proof Trivial.\fp

\begin{theorem}[Level 2]
  \begin{enumerate}
  \item The functions
  $$f_{2,j}(\tau)=4E_{1,-4}(\tau)+(-1)^j2\sqrt{-1}(E_{1,-4}(\tau/2)-E_{1,-4}((\tau+1)/2))$$
(where we recall that $E_{1,-4}=\th^2/4$) have the following properties:
  \begin{align*}
    &f_{2,j}^4(\tau)=F_{2,j}(\tau)\;,\\
    &f_{2,j}(\tau)f_{2,j+1}(\tau)=E_{2,2}(\tau)\;,\quad f_{2,j}(((-1)^j+\sqrt{-1})/2)=0\;,\\
    &f_{2,j}(\tau+1)=f_{2,j+1}(\tau)\;,\quad f_{2,j}(-1/(2\tau))=((-1)^j-\sqrt{-1})\tau f_{2,j+1}(\tau)\;.\end{align*}
\item The functions
  $$g_{2,j}(\tau)=\sum_{0\le i\le 3}(1-2\delta_{i,j})E_{1,-8}((\tau+i)/4)$$
  have the following properties:
  \begin{align*}
    &g_{2,j}^2(\tau)=G_{2,j}(\tau)\;,\\
    &\prod_{0\le j\le 3}g_{2,j}(\tau)=E_{2,4}(\tau)\;,\quad g_{2,j}((4-j)+\sqrt{-2}/2)=0\,\\
    &g_{2,j}(\tau+1)=g_{2,j+1}(\tau)\;,\quad g_{2,j}(-1/(2\tau))=\eps_j\sqrt{-2}\tau g_{2,j}(\tau)\;,
  \end{align*}
  with $\eps_j=1$ if $4\mid j$ and $\eps_j=-1$ if $4\nmid j$.
  \end{enumerate}
\end{theorem}

\Proof As in level 1, one can do a non-illuminating proof by working in
finite-dimensional vector spaces of modular forms for $\G(8)$.\fp

This theorem shows both that the functions
$E_{2,4}(\tau)-(-1)^j\sqrt{-256}\eta_2^4(\tau)$ are
fourth powers of natural modular forms, and that
$E_{2,2}(\tau)=2E_2(2\tau)-E_2(\tau)$ is the product of two conjugate
Eisenstein series of weight $1$.

\section{Results in Level $3$}

For simplicity, we set $\eta_3(\tau)=\eta(\tau)\eta(3\tau)$, of weight $1$
(not to be confused with $\eta(3\tau)$ itself).

In this level, first note the identity $E_{3,2}=(6E_{1,-3})^2$. In addition,
note that the function $E_{3,4}^2/E_{3,2}$ is a modular form (i.e., is
holomorphic) since one easily shows that
$E_{3,4}^2/E_{3,2}=E_{3,2}^3-108\eta_3^6$. Even better, its square root
$E_{3,4}/E_{3,2}^{1/2}\in M_3(\G_0(3),\chi_{-3})$, since
$E_{3,4}/E_{3,2}^{1/2}=-9(E_{3,-3,1}+3E_{3,1,-3})$, where
  $$E_{3,-3,1}=-1/9+\sum_{n\ge1}\bigl(\sum_{d\mid n}d^2\lgs{-3}{d}\bigr)q^n\text{\quad and\quad}
  E_{3,1,-3}=\sum_{n\ge1}\bigl(\sum_{d\mid n}d^2\lgs{-3}{n/d}\bigr)q^n$$
  are the two normalized Eisenstein series in $M_3(\G_0(3),\chi_{-3})$.

  Thus, we have two basic formulas in level $3$:
$$E_{3,2}^3-E_{3,4}^2/E_{3,2}=108\eta_3^6\text{\quad and\quad}E_{3,2}^4-E_{3,4}^2=3888\eta_3^6E_{1,-3}^2\;,$$
  each of which can be factored in two ways and leads to the following:
\begin{definition}
  We define
  \begin{align*}
    F_{3,j}(\tau)&=(E_{3,4}/E_{3,2}^{1/2})(\tau)-(-1)^j6\sqrt{-3}\eta_3^3(\tau)\;,\\
    G_{3,j}(\tau)&=E_{3,2}(\tau)-3\root3\of{4}\rho^j\eta_3^2(\tau)\;,\\
    F'_{3,j}(\tau)&=E_{3,4}(\tau)-(-1)^j36\sqrt{-3}(\eta_3^3E_{1,-3})(\tau)\;,\\
    G'_{3,j}(\tau)&=E^2_{3,2}(\tau)-(-1)^j36\sqrt{3}(\eta_3^3E_{1,-3})(\tau)\;.\end{align*}
\end{definition}

We have evidently $F'_{3,j}=6E_{1,3}F_{3,j}$.

\begin{proposition}
  We have
  \begin{align*}
    &F_{3,j}(\tau+1)=F_{3,j+1}(\tau)\;,\quad F_{3,j}(-1/(3\tau))=-3\sqrt{-3}\tau^3F_{3,j+1}(\tau)\;,\\
    &F_{3,j}(\tau)F_{3,j+1}(\tau)=E_{3,2}^3(\tau)\;,\\
    &G_{3,j}(\tau+1)=G_{3,j+1}(\tau)\;,\quad G_{3,j}(-1/(3\tau))=-3\tau^2G_{3,j}(\tau)\;,\\
    &G_{3,j}(\tau)G_{3,j+1}(\tau)G_{3,j+2}(\tau)=(E_{3,4}^2/E_{3,2})(\tau)\;,\\
    &F'_{3,j}(\tau+1)=F'_{3,j+1}(\tau)\;,\quad F'_{3,j}(-1/(3\tau))=-9\tau^4F'_{3,j+1}(\tau)\;,\\
    &F'_{3,j}(\tau)F'_{3,j+1}(\tau)=E_{3,2}^4(\tau)\;,\\
    &G'_{3,j}(\tau+1)=G'_{3,j+1}(\tau)\;,\quad G'_{3,j}(-1/(3\tau))=9\tau^4G'_{3,j}(\tau)\;,\\
    &G'_{3,j}(\tau)G'_{3,j+1}(\tau)=E_{3,4}^2(\tau)\;.\end{align*}
\end{proposition}

\Proof Trivial.\fp

\begin{theorem}[Level 3]
  The functions
  $$f_{3,j}(\tau)=6E_{1,-3}(2\tau)+(-1)^j\sqrt{-3}(E_{1,-3}(\tau/2)-E_{1,-3}((\tau+1)/2))$$
  have the following properties:
  \begin{align*}
    &f_{3,j}^3(\tau)=F_{3,j}(\tau)\;,\\
    &f_{3,j}(\tau)f_{3,j+1}(\tau)=E_{3,2}(\tau)\;,\quad f_{3,j}(((-1)^j3+\sqrt{-3})/6)=0\;,\\
    &f_{3,j}(\tau+1)=f_{3,j+1}(\tau)\;,\quad f_{3,j}(-1/(3\tau))=(((-1)^j3-\sqrt{-3})/2)\tau f_{3,j+1}(\tau)\;.\end{align*}
  \end{theorem}

\Proof As in previous levels, work over $\G(6)$.\fp

This shows both that $F_{3,j}$ is the cube of a modular form, and that
$(3E_2(3\tau)-E_2(\tau))/2$ is the product of two conjugate weight $1$ Eisenstein series.

Contrary to levels $1$, $2$, and $4$, I have not been able to find
corresponding results for the functions $G_{3,j}$, $F'_{3,j}$, and $G'_{3,j}$.

\section{Results in Level $4$}

For simplicity, we set $\eta_4(\tau)=\eta(\tau)\eta(4\tau)/\eta(2\tau)$, of
weight $1/2$ (not to be confused with $\eta(4\tau)$ itself). The space of
Eisenstein series of weight $2$ is now two-dimensional (it was one-dimensional
in levels $2$ and $3$), so in addition to $E_{4,2}$ we need to introduce
an additional Eisenstein series, and we choose
$$E'_{4,2}(\tau)=4E_2(4\tau)-4E_2(2\tau)+E_2(\tau)$$
(we will not use derivatives in this paper so this notation will not lead to
any confusion).

This level has many more eta quotient identities:
\begin{align*}
  E_{4,2}(\tau)=\th^4(\tau)&=\eta(2\tau)^{20}/(\eta(\tau)\eta(4\tau))^8\;,\\
  E_{4,2}(\tau)-E'_{4,2}(\tau)&=32\eta(4\tau)^8/\eta(2\tau)^4\;,\\
  E_{4,2}(\tau)+E'_{4,2}(\tau)&=2\eta(\tau)^8/\eta(2\tau)^4\;,\\
  E_{4,2}^2(\tau)-{E'_{4,2}}^2(\tau)&=64(\eta(\tau)\eta(4\tau)/\eta(2\tau))^8=64\eta_4^8(\tau)\;,\\
  \eta(2\tau)^{24}&=(\eta(\tau)\eta(4\tau))^8(\eta(\tau)^8+16\eta(4\tau)^8)\;.
  \end{align*}

The basic identity $E^2_{4,2}-{E'_{4,2}}^2=64\eta_4^8(\tau)$ can thus be
factored in two additional ways, and leads to the following:
\begin{definition}
  We define
  $$F_{4,j}(\tau)=E'_{4,2}(\tau)-(-1)^j\sqrt{-64}\eta_4^4(\tau)\text{\quad and\quad}G_{4,j}(\tau)=E_{4,2}(\tau)-(-1)^j8\eta_4^4(\tau)\;.$$
\end{definition}

\begin{proposition}
  We have
  \begin{align*}
  &F_{4,j}(\tau+1)=F_{4,j+1}(\tau)\;,\quad F_{4,j}(-1/(4\tau))=4\tau^2F_{4,j+1}(\tau)\;,\\
  &F_{4,j}(\tau)F_{4,j+1}(\tau)=E_{4,2}^2(\tau)\;,\\
  &G_{4,j}(\tau+1)=G_{4,j+1}(\tau)\;,\quad G_{4,j}(-1/(4\tau))=-4\tau^2G_{4,j}(\tau)\;,\\
  &G_{4,j}(\tau)G_{4,j+1}(\tau)={E'_{4,2}}^2(\tau)\;.\end{align*}
\end{proposition}

\Proof Trivial.\fp

\begin{theorem}[Level 4]
  \begin{enumerate}\item
  The functions $f_{4,j}(\tau)=\th(\tau/2-(-1)^j/4)$
  have the following properties:
  \begin{align*}
    &f_{4,j}^4(\tau)=F_{4,j}(\tau)\;,\\
    &f_{4,j}(\tau)f_{4,j+1}(\tau)=E_{4,2}^{1/2}(\tau)=\th^2(\tau)\;,\\
    &f_{4,j}(\tau+1)=f_{4,j+1}(\tau)\;,\quad f_{4,j}(-1/(4\tau))=(1+(-1)^j\sqrt{-1})(\tau/i)^{1/2} f_{4,j+1}(\tau)\;.\end{align*}
\item
  The functions $g_{4,j}(\tau)=2\th^2(2\tau)-\th^2((\tau+j)/2)$
  have the following properties:
  \begin{align*}
    &g_{4,j}^2(\tau)=G_{4,j}(\tau)\;,\\
    &g_{4,j}(\tau)g_{4,j+1}(\tau)=E'_{4,2}(\tau)\;,\\
    &g_{4,j}(\tau+1)=g_{4,j+1}(\tau)\;,\quad g_{4,j}(-1/(4\tau))=(-1)^j2\sqrt{-1}\tau g_{4,j}(\tau)\;.\end{align*}
  \end{enumerate}
  \end{theorem}

\Proof As usual, work over $\G(8)$.\fp

This shows both that $F_{4,j}$ is the fourth power of a natural function,
and that $((4E_2(4\tau)-E_2(\tau))/3)^{1/2}$ is a product of two conjugate
unary theta series of weight $1/2$. Perhaps even more surprising is the
equality $f_{4,j}f_{4,j+1}=\th^2$, which shows that this function is a product
of two forms of weight $1/2$ in two different ways. Note that this identity
also follows trivially from the product formula
$\th(\tau)=\prod_{n\ge1}(1-q^{2n})(1+q^{2n-1})^2$ by rearranging terms.

\section{Conclusion}

The existence of the above identities leads to several questions.
First, it is quite plausible that similar identities exist in higher levels,
for instance in levels $N$ such that $(N-1)\mid 24$, or more generally such
that $\G_0^*(N)$ has genus $0$.

Second, as suggested by D.~Zagier, are these factorizations unique, more
generally to what extent do holomorphic modular forms have unique
factorizations (up to scalars) into irreducible ones, is there a description
of the \emph{complete} set of holomorphic modular forms dividing a given one,
for instance $E_4$ ?

\end{document}